\documentstyle{amsppt}

\NoBlackBoxes

\topmatter
\title
On  Differences of Riesz Homomorphisms
\endtitle
\author    S.~S.~Kutateladze\endauthor
\thanks  A talk at the United Mathematical Seminar of 
the Sobolev Institute
on February 4, 2005
\endthanks
\endtopmatter

\document
\noindent
Which closed hyperplanes in a Banach lattice  are closed
Riesz subspaces, i.e. closed under finite meets and joins?
It turns out that these are exactly the kernels of
differences of Riesz homomorphisms on the initial Riesz space.
The motivation behind this question
is the celebrated Stone Theorem about the structure of
Riesz subspaces in the Banach lattice
$C(Q,\Bbb R)$
of continuous real functions on a compact space ~$Q$.
This theorem may be rephrased in the terms we are interested
in as follows:

\proclaim{Stone Theorem}
Each closed Riesz subspace of $C(Q,\Bbb R)$
is the intersection of the kernels
of some differences of  Riesz homomorphisms on~$C(Q,\Bbb R)$.
\endproclaim

In view of this theorem we find it reasonable to
refer to a difference of Riesz homomorphisms
on a Riesz space ~$X$ as a~
{\it two-point relation\/} on~$X$.
We are not obliged to assume here that
the Riesz homomorphisms under study act into the reals~$\Bbb R$.
In what follows we will consider the general
Riesz homomorphisms from ~$X$ to an arbitrary
Kantorovich (= Dedekind complete Riesz) space with
base a~ complete Boolean algebra~$B$.

We are not obliged to assume here that
the Riesz homomorphisms under study act into the reals~$\Bbb R$.
In what follows we will consider the general
Riesz homomorphisms from ~$X$ to an arbitrary
Kantorovich (= Dedekind complete Riesz) space $Y$ with
base a~ complete Boolean algebra~$B$.

So, let $T:X\to Y$ be an order bounded operator
whose positive and negative parts are Riesz homomorphisms.
Observe that for every  band projection
$b\in B$ the operator $bT$, called a~{\it stratum\/} of~$T$ is a Riesz subspace of~$X$.
In fact, the converse is valid too. In other words, we have the following

\proclaim{Main Theorem}
An order bounded operator from a Riesz space to a Kantorovich space
is a two-point relation if and only if the kernel of its every
stratum is a Riesz subspace
of the ambient Riesz space.
\endproclaim

The analysis below is rather transparent if we
apply a well-developed technique of
``nonstandard scalarization.'' This technique
implements the Kantorovich heuristic principle
and reduces operator problems to the case of functionals.

The sketch of the proof is as follows:
Using the functors of canonical embedding and descent to the
Boolean valued universe $\Bbb V^{(B)}$,
we reduce the matter to characterizing
scalar two-point relations on Riesz spaces
over dense subrings of the reals~$\Bbb R$.
To solve the resulting scalar problem, we use
one of the formulas of subdifferential calculus
which is known as

\proclaim{Decomposition Theorem}
Assume that $H_1,\dots,H_N$ are cones in a Riesz space~$X$.
Assume further that $f$ and $g$ are positive functionals on ~$X$.

The inequality
$$
f(h_1\vee\dots\vee h_N)\ge g(h_1\vee\dots\vee h_N)
$$
holds for all
$h_k\in H_k$ $(k:=1,\dots,N)$
if and only if to each decomposition
of~ $g$ into a sum of~$N$ positive terms
$g=g_1+\dots+g_N$
there is a decomposition of ~$f$ into a sum of~$N$
positive terms $f=f_1+\dots+f_N$
such that
$$
f_k(h_k)\ge g_k(h_k)\quad
(h_k\in H_k;\ k:=1,\dots,N).
$$
\endproclaim

Obviously, the Riesz space  in this
theorem may be viewed over an arbitrary dense subfield
~$R$ of the reals~$\Bbb R$.
We know in passing that the results to follow
are preserved in the general modules admitting convex analysis.

Descending the Main Theorem from an appropriate
Boolean valued universe
or, which is equivalent, using the characterization of the modules
admitting convex analysis,
we can arrive to an analogous description for
a dominated module homomorphism with
kernel a Riesz subspace in modules over
an almost rational subring of the orthomorphism ring of the
range.

From the Main Theorem it is immediate that the Stone Theorem cannot
be abstracted beyond the limits of $AM$-spaces.
Indeed, if each closed Riesz subspace of a Banach lattice is
an intersection of two-point relations then there are sufficiently
many Riesz homomorphisms to separate the points of the
Banach lattice under consideration.

\enddocument